# Some mathematical properties of Ch. Bartlett's 'chi ratio'

**Definition of chi**

Based on artistic interpretations, art professor Christopher Bartlett (Towson University, USA) suggested introducing a new mathematical constant related to the golden ratio (see [1]). He called it the 'chi ratio', as χ follows ϕ in the Greek alphabet. Without going too much into the artistic considerations, which may be subject to criticism similar to the 'golden ratio debunking', we focus on showing that the chi ratio is interesting as a number as such, with remarkable geometric properties, just as the golden ratio is.

Indeed, one of the reasons why the golden section ϕ = (1+√5)/2 = 1.618… is a pleasing mathematical number, is that a rectangle of width 1 and length ϕ can subdivided into a square and another rectangle proportional to the original rectangle. That smaller proportional rectangle will have width 1/ϕ and length 1, and the square a side 1 because ϕ = 1/ϕ + 1 (see figure 1). They can be easily constructed using a diagonal and the perpendicular from a vertex not on that diagonal. Of course, there are (many) other methods to get ϕ, but the latter will be applied here in other instances.

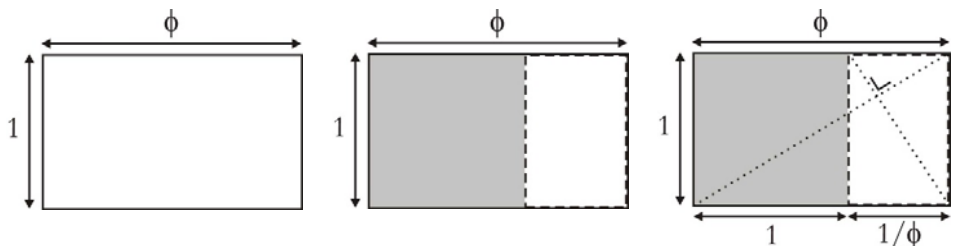

**Figure 1**: A rectangle with proportion ϕ can be subdivided in a square (grey) and a rectangle proportional to the original rectangle (dashed line).

Bartlett proposed to continue this procedure, that is, he wondered what kind of rectangle would be constructed in a similar way, if the starting point were a golden rectangle, instead of a square? He proposed to consider a rectangle with width 1 and length ϕ or one with width 1/ϕ and length 1, in order to extend the 'golden rectangle' by a rectangle of the same proportions as the combination of it with a golden section rectangle. Again, the construction can be easily checked, using a diagonal and the perpendicular from a vertex not on that diagonal. In the case where the original rectangle has width 1 and length ϕ, Bartlett used the notation 1/χ' for the width of the added rectangle (while it will have length 1). Thus, the combination of both will have width 1 and length χ'. In the case where the original rectangle has width 1/ϕ and length 1, he denoted the width of the added rectangle by 1/χ (again with length 1). Now, the combination of both will have width 1 and length χ (see figure 2).

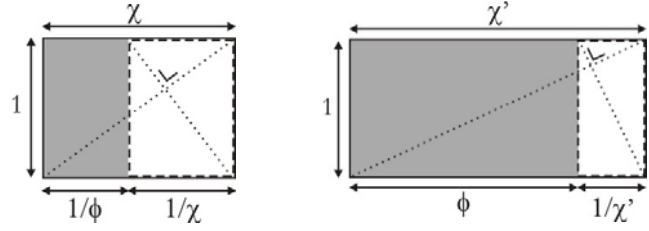

**Figure 2**: Subdividing a rectangle with proportions χ or χ' in a golden rectangle (grey) and a rectangle proportional to the original rectangle (dashed line).

The value of χ follows from the equation $\chi = 1/\phi + 1/\chi$, as it transforms in the quadratic equation $\chi^2 - (1/\phi)\chi - 1 = 0$. Thus, $\chi = \frac{(1/\phi) \pm \sqrt{(1/\phi)^2 + 4}}{2} = \frac{\sqrt{5} - 1 \pm \sqrt{22 - 2\sqrt{5}}}{4}$. The positive solution is 1.355… and this corresponds to what Bartlett called 'the chi ratio': χ =1.355… .

Similarly, χ' follows from the equation $\chi' = \phi + 1/\chi'$, transforming into $\chi'^2 - \phi\chi' - 1 = 0$, so that $\chi' = \frac{\phi \pm \sqrt{(1/\phi)^2 + 4}}{2} = \frac{1 + \sqrt{5} \pm \sqrt{22 + 2\sqrt{5}}}{4}$. Now, the positive solution is 2.095…. and Bartlett denoted it by χ' = 2.095…

In the last century, and only in the last century, the golden section has gained an enormous popularity in (pseudo-) science, due to the works of the German psychologist Adolf Zeising (1810-1876) and the Romanian diplomat Prince Matila Ghyka (1881 – 1965). They imagined the golden section was used for esthetic purposes as well, that is, that the proportion 1.618… would have been omnipresent in art and architecture. There are no historic or scientific arguments for it (see [4]), but perhaps the use of Bartlett's chi ratio in art is not that odd. Dutch architect Dom Hans van der Laan's 'plastic number' ψ = 1.324…. (1928; see [7]) and Spanish architect Rafael de la Hoz's 'Cordovan proportion' $c$ = 1.306… (1973; see [3]) are both close to the classical proportion of 4/3, called the 'sesquitertia' (see [6]). Unquestionably, the 'sesquitertia' was used by Vitruvius, Pacioli and Leonardo, in contrast to the golden section. Moreover, since Bartlett (unconsciously!) followed the right angle construction that was so treasured by Le Corbusier, his discovery might resist future 'chi ratio debunkings' (see [2]). Still, the goal of the present paper is the mathematical properties of χ, and so we will limit ourselves to those.

**Some immediate properties**

From

$$\chi^2 - (1/\phi)\chi - 1 = 0 \qquad \text{and} \qquad \chi'^2 - \phi\chi' - 1 = 0$$

follows that

$$\chi = \sqrt{1 + (1/\phi)\chi} \qquad \text{and} \qquad \chi' = \sqrt{1 + \phi\chi'}$$

Thus

$$\chi = \sqrt{1 + \phi^{-1}\sqrt{1 + \phi^{-1}\sqrt{1 + \ldots}}} \quad \text{and} \quad \chi' = \sqrt{1 + \phi\sqrt{1 + \phi\sqrt{1 + \cdots}}}$$

This can be compared to the well-known root expression for $\phi$:

$$\phi = \sqrt{1 + \sqrt{1 + \sqrt{1 + \cdots}}}.$$

Also, from

$$\chi^2 - (1/\phi)\chi - 1 = 0 \quad \text{and} \quad \chi'^2 - \phi\chi' - 1 = 0$$

it follows that

$$\chi = (1/\phi) + 1/\chi \quad \text{and} \quad \chi' = \phi + 1/\chi'$$

Thus

$$\chi = \left(\tfrac{1}{\phi}\right) + \cfrac{1}{\left(\tfrac{1}{\phi}\right) + \cfrac{1}{\left(\tfrac{1}{\phi}\right) + \cdots}} \quad \text{and} \quad \chi' = \phi + \cfrac{1}{\phi + \cfrac{1}{\phi + \cdots}}$$

This can be compared to the well-known continued fraction for $\phi$:

$$\phi = 1 + \cfrac{1}{1 + \cfrac{1}{1 + \cdots}}.$$

A geometric interpretation is given in figure 3.

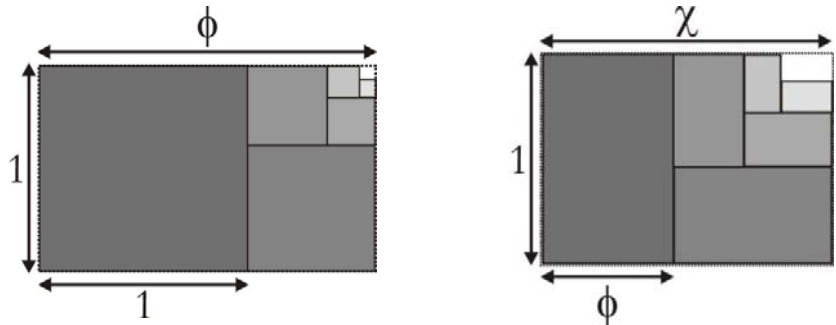

**Figure 3**: Interpretations of the continued fractions for $\varphi$ and for $\chi$.

The fractions 2, 3/2, 5/3, 8/5, … successively taken from the latter continued fraction give rise to the Fibonacci sequence ($F_n$): (1, 1,) 2, 3, 5, 8, … in which each term is the sum of the previous ones. We follow a similar approach with the continued fraction for $\chi'$ and consider the fraction:

$$\phi+1, \; \phi + \cfrac{1}{\phi+1}, \; \phi + \cfrac{1}{\phi + \cfrac{1}{\phi+1}}, \; \phi + \cfrac{1}{\phi + \cfrac{1}{\phi + \cfrac{1}{\phi+1}}}, \; \phi + \cfrac{1}{\phi + \cfrac{1}{\phi + \cfrac{1}{\phi + \cfrac{1}{\phi+1}}}}, \; \ldots$$

That is, taking into account that $\phi^2 = \phi+1$:

$\phi+1,$

$\phi + \cfrac{1}{\phi+1} = (\phi^2+\phi+1)/(\phi+1) = (2\phi+2)/(\phi+1) = 2,$

$$\phi + \cfrac{1}{\phi+\cfrac{1}{\phi+1}} = \phi+1/(2) = (2\phi+1)/2,$$

$$\phi + \cfrac{1}{\phi+\cfrac{1}{\phi+\cfrac{1}{\phi+1}}} = \phi+2/(2\phi+1) = (2\phi^2+\phi+2)/(2\phi+1) = (3\phi+4)/(2\phi+1),$$

$$\phi + \cfrac{1}{\phi+\cfrac{1}{\phi+\cfrac{1}{\phi+\cfrac{1}{\phi+1}}}} = \phi+(2\phi+1)/(3\phi+4) = ((3\phi^2+4\phi)+(2\phi+1))/(3\phi+4)$$

$$= ((3\phi+3+(6\phi+1))/(3\phi+4) = (9\phi+4)/(3\phi+4)$$

…

This gives rise to a series $(H_n)$: $(\phi+1, 2,)$ $2\phi+1$, $3\phi+4$, $9\phi+4$, $16\phi+13$, $38\phi+20$, $74\phi+51$,… and from its construction it follows the series $H_{n+1}/H_n$ converges to $\chi' = 2.095…$ Its rule is simple and somehow similar to the Fibonacci sequence: if $a\phi+b$, $c\phi+d$ are consecutive terms, than the next one is $(a+c+d)\phi+(b+c)$ (Note that if the more obvious rule, $(a\phi+b)+(c\phi+d)$ as in the Fibonacci series, is used, the series $H_{n+1}/H_n$ would converge to $\phi$). It is not an integer sequence, in contrast to the Fibonacci sequence. The first numbers in the series are: $(2.618…, 2,)$ $4.236…$, $8.854…$, $18.562…$, $38.888…$, $81.484…$, $170.732…$, … .

As $\phi = 1/2 + \sqrt{5}/2$ and $\chi = \sqrt{5}/4 - 1/4 + (\sqrt{22-2\sqrt{5}})/4$ are rational expressions with square roots, they can be constructed with compass and ruler in a finite number of steps (see figure 4).

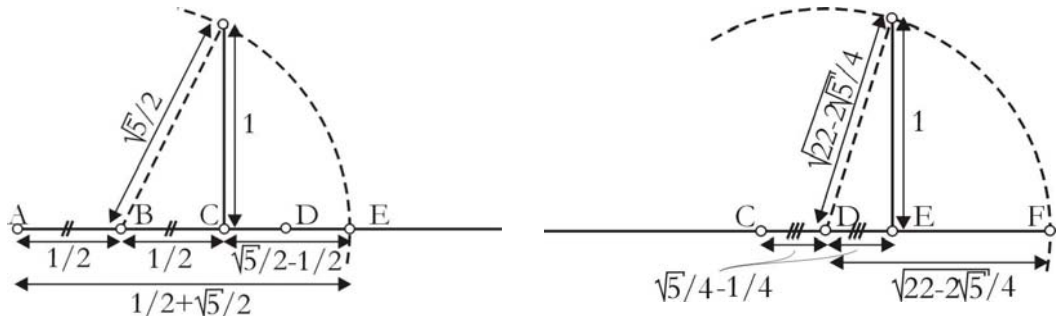

**Figure 4**: Construction of the $\varphi$ - and $\chi$-ratio.

From the above equations it follows that the solutions of $x^2 - (1/\phi)x - 1 = 0$ are $\chi$ and $-1/\chi$, while the solutions of $x^2 - \phi x - 1 = 0$ are $\chi'$ and $-1/\chi'$, and thus:

$$(x-\chi)(x+\frac{1}{\chi})(x-\chi')\left(x+\frac{1}{\chi'}\right) = \left(x^2 - (\chi-\frac{1}{\chi})\cdot x - 1\right)\left(x^2 - (\chi'-\frac{1}{\chi'})\cdot x - 1\right)$$

$= (x^2 - \phi x - 1)\cdot(x^2 - (1/\phi)x - 1)$
$= x^4 - (\phi +1/\phi)x^3 - x^2 + (\phi +1/\phi)x + 1$

However, in order to find the so-called 'minimal polynomial' of $\chi$ and get rational coefficient for $x$, a change of signs provides the solution (this was noticed by A. Redondo Buitrago [6], while discussing this paper):

$$(x-\chi)(x+\frac{1}{\chi})(x+\chi')\left(x-\frac{1}{\chi'}\right) = \left(x^2 - (\chi - \frac{1}{\chi})\cdot x - 1\right)\left(x^2 + (\chi' - \frac{1}{\chi'})\cdot x - 1\right)$$

$$= (x^2 - \phi x - 1)\cdot(x^2 + (1/\phi)x - 1)$$

$$= x^4 - (\phi - 1/\phi)x - 3x^2 + (\phi - 1/\phi)x + 1$$
$$= x^4 - x - 3x^2 + x + 1$$

Thus, $\phi$ and $\chi$ combine in a pleasant way from a purely mathematical point of view.

**Extending arbitrary rectangles proportionally**

Bartlett's method with the orthogonal diagonals linking a square to a golden section rectangle, and a golden section rectangle to a chi ratio rectangle of course works equally well when starting with an arbitrary rectangle. In case $x < \phi$, the proportion of the longer side to the shorter side gives rise to the equation $1/(x - 1/x) = x/(x^2 - 1) = \rho$, if $\rho$ is the given proportion. In case $x > \phi$, the proportion of the longer side to the shorter side gives rise to the equation $(x - 1/x)/1 = (x^2 - 1)/x = \rho$, if $\rho$ is the given proportion (see figure 5).

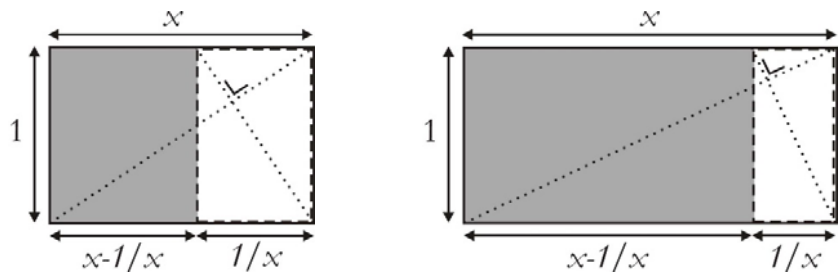

**Figure 5**: Finding proportional rectangles in case $x < \phi$ (left) and $x > \phi$ (right).

For instance, if one starts with a 1 by 2 rectangle (1 and 2 are the 2$^{nd}$ and 3$^{th}$ Fibonacci numbers), $\rho = 2/1$, and thus an $x = 1+\sqrt{2}$ rectangle will provide the desired subdivision of a 1 by $(1+\sqrt{2})$ rectangle in a 1 by 2 rectangle and a 1 by $1/(1+\sqrt{2})$ rectangle (see figure 6). The number $1+\sqrt{2}$ is called the 'silver section'.

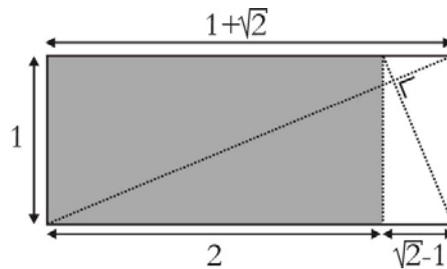

**Figure 6**: A 1 by 2 rectangle extended by and as a *silver* rectangle.

If one starts with a 2 by 3 rectangle (2 and 3 are the 3$^{th}$ and 4$^{th}$ Fibonacci numbers), $\rho = 3/2$, and thus an $x = 2$ rectangle will provide the desired subdivision in a 1 by 3/2 rectangle and a 1 by 1/2 rectangle. This reminds one the often occurring subdivisions of (parts of) a canvas in 4 parts, as in

Georges Seurat's painting 'Parade de cirque' (see figure 7). Indeed, following Marguerite Neveux (France), there is no trace of a golden section interpretation in this painting, but, on the contrary, of a 2-4-8 division. Golden section critic Neveux asserted Seurat subdivided the left part of the painting in 4 parts (see [5]), and this fits well in Bartlett's right angle subdivision method.

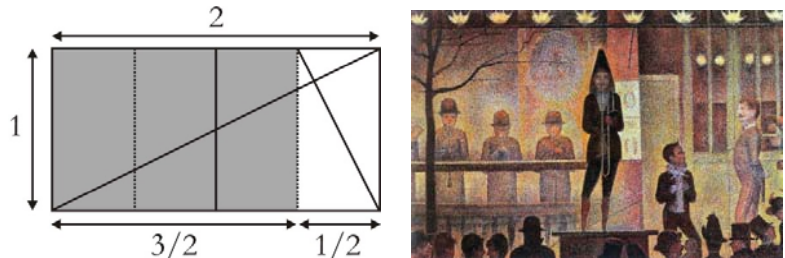

**Figure 7**: Proportional rectangles can explain Seurat's painting 'Parade de cirque'.

We can generalize Bartlett's method: what about extending a chi ratio rectangle with a rectangle proportional to that rectangle together with the chi ratio rectangle? It yields a rectangle with length 1.434… (and width 1). And what if this method is continued 'at infinity', until the equation $x/(x^2 -1) = x$ is obtained? Surprisingly – or not – the positive solution is the square root of 2 rectangle, that is, the (European) DIN A4 paper format (see figure 8).

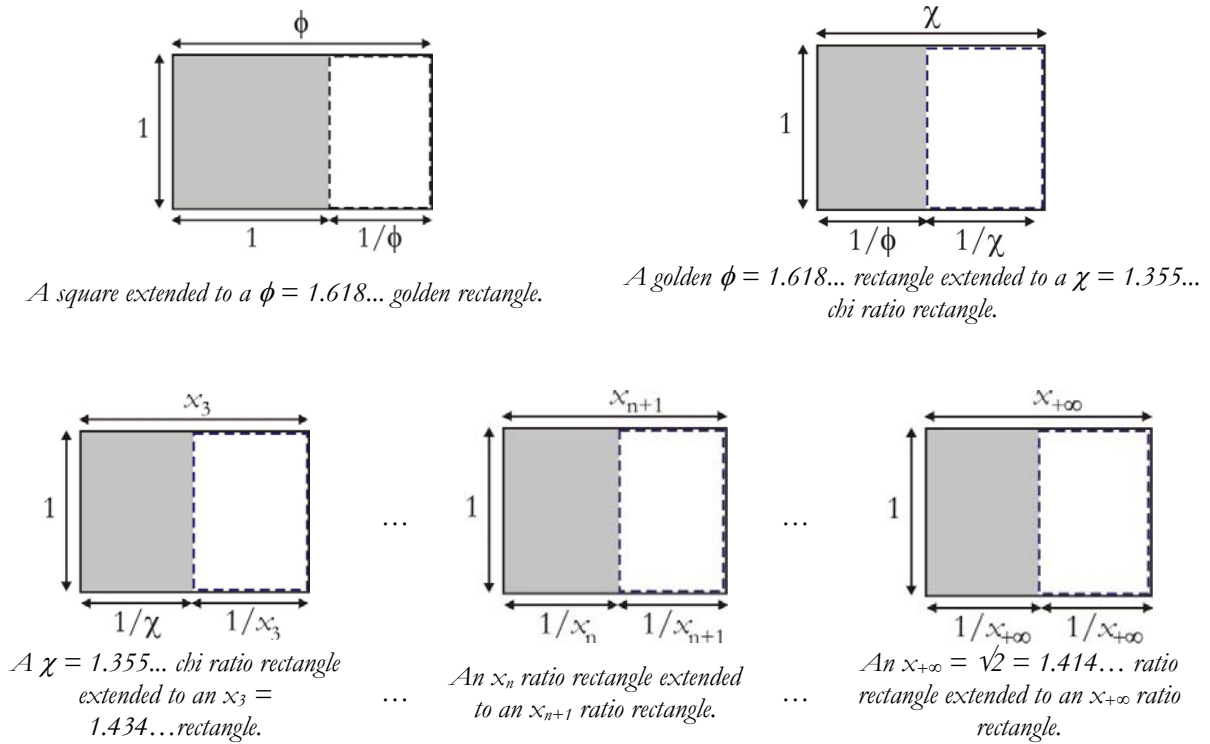

*A square extended to a $\phi$ = 1.618... golden rectangle.*

*A golden $\phi$ = 1.618... rectangle extended to a $\chi$ = 1.355... chi ratio rectangle.*

*A $\chi$ = 1.355... chi ratio rectangle extended to an $x_3$ = 1.434…rectangle.*

*An $x_n$ ratio rectangle extended to an $x_{n+1}$ ratio rectangle.*

*An $x_{+\infty}$ = $\sqrt{2}$ = 1.414… ratio rectangle extended to an $x_{+\infty}$ ratio rectangle.*

**Figure 8**: Extending Bartlett's method at infinitum.

**Paper folding**

A surprising advantage of Bartlett's diagonal construction to create proportional rectangles is that the method can be executed easily by folding, without any instrument or calculator. Say a sheet of paper has width 1 and an arbitrary length $x$. After folding one diagonal, the perpendicular through one vertex not on that diagonal can be obtained easily by folding the paper such that the folds correspond (see figure 9). It yields a length $1/x$.

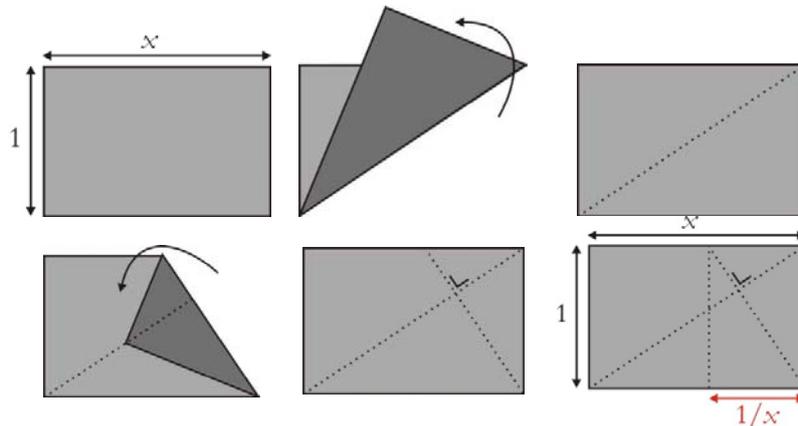

**Figure 9**: Folding a proportional rectangle, given an arbitrary rectangle.

As a square can be easily folded (see figure 10), a sequence of squares gives rectangles of integer length (given the width is 1). Repetition of the folding method to a rectangle of width 1 and length $n$, gives a rectangle of width 1 and length $n+1/(n+1/n)$, or as many steps as desired in the expansion of the continued fraction $n + \cfrac{1}{n+\cfrac{1}{n+\cfrac{1}{n+\dots}}}$. For instance, if $n = 2$, $2 + \cfrac{1}{2+\cfrac{1}{2+\cfrac{1}{2+\dots}}}$ will correspond to the silver section $1 + \sqrt{2}$, and thus approximations for that irrational number based on its continued fraction can be folded easily. The folding process based on $n=3$ given in the illustration will converge to the so-called 'bronze mean', $(3 + \sqrt{13})/2$.

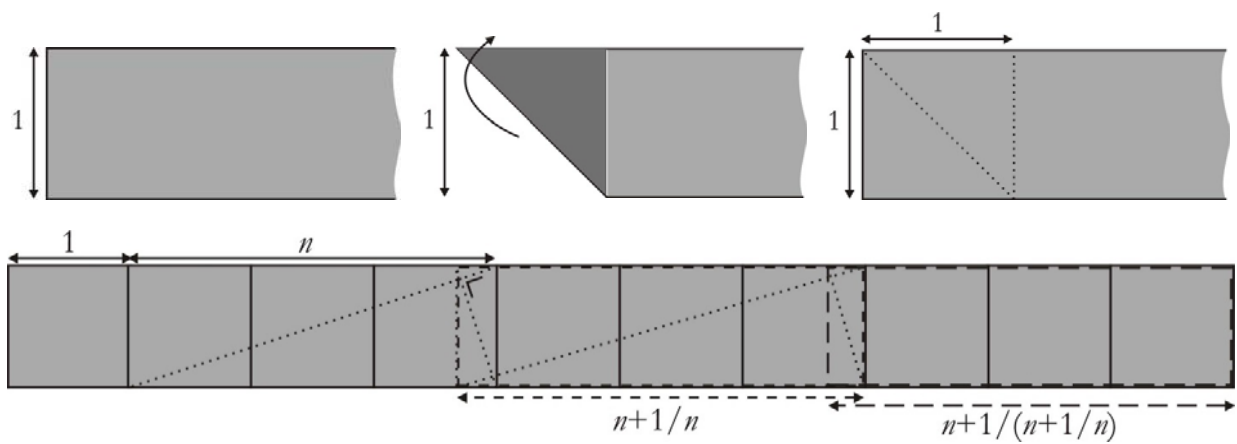

**Figure 10**: Folding a square is easy, and thus so is a rectangle of any integer length; hence a continued fraction can be folded to a desired degree of precision (here, obviously, $n = 3$).

An 'application' is the construction of harmonic mean $H = 2/(1/m + 1/n)$ given $m$ and $n$. Indeed, $(1/m + 1/n) = (m+n)/mn$ easily follows from the above constructions, and thus so does $H$ (see figure 11).

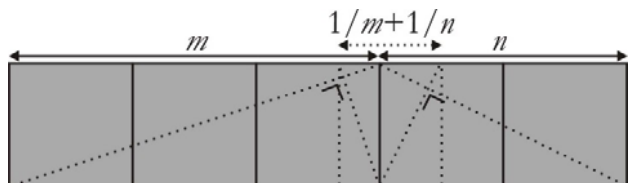

**Figure 11**: Folding $(m+n)/mn$ given $m$ and $n$ (in the figure, $m=3$ and $n=2$).

The folding procedure works for a rectangle with an arbitrary length. Moreover, repetition of the folding method starting from a rectangle with width 1 and arbitrary length $x$ will, 'in the end', provide the golden section (see figure 12). And thus, perhaps there are some good elements in the 'golden section myth' after all.

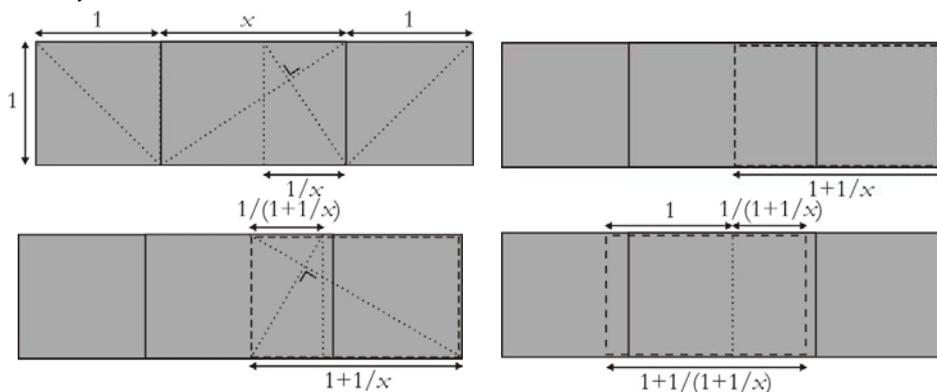

**Figure 12**: Folding a rectangle with width $x$ yields a rectangle with width $1+1/(1+(1/1… 1+/(1+x)…)$, an approximation for ϕ.